% X-sliced-and-diced-by: 'savemail' 0.1, Aug 30, 1997

% \def\cal{\Cal}

\def\bf{\fam \bffam \tenbf }
\input amssym

%%%\magnification=\magstep1
\hsize=15,4truecm
\vsize=23.5truecm
\mathsurround=1pt

\def\chapter#1{\par\bigbreak \centerline{\bf #1}\medskip}

\def\section#1{\par\bigbreak {\bf #1}\nobreak\enspace}

\def\sqr#1#2{{\vcenter{\hrule height.#2pt
      \hbox{\vrule width.#2pt height#1pt \kern#1pt
         \vrule width.#2pt}
       \hrule height.#2pt}}}

\def\k{\kappa}
\def\o{\omega}

\def\d{\delta}

\def\l{\lambda}

\def\a{\alpha}
\def\b{\beta}

\def\g{\gamma}

\def\n{\eta}

%=====================

\def\A{{\cal A}}
\def\B{{\cal B}}
\def\C{{\cal C}}

\def\M{{\bf M}}

%========================

%\def\pm{\buildrel{\scriptscriptstyle +}\over {\scriptscriptstyle -}}
\def\th #1 #2. #3\par\par{\medbreak{\bf#1 #2.
\enspace}{\sl#3\par}\par\medbreak}
\def\rem #1 #2. #3\par{\medbreak{\bf #1 #2.
\enspace}{#3}\par\medbreak}
\def\proof{{\bf Proof}.\enspace}
\def\sqr#1#2{{\vcenter{\hrule height.#2pt
      \hbox{\vrule width.#2pt height#1pt \kern#1pt
         \vrule width.#2pt}
       \hrule height.#2pt}}}
\def\eop{\mathchoice\sqr34\sqr34\sqr{2.1}3\sqr{1.5}3}

%====================================================================%
                                                                     %
% A macro for making references and blocks.                          %
                                                                     %
\newdimen\refindent\newdimen\plusindent                              %
\newdimen\refskip\newdimen\tempindent                                %
\newdimen\extraindent                                                %
                                                                     %
%\refskip has to be determined by the user! Otherwise \parindent is  %
%used, in accordance with \item.                                     %
                                                                     %
\def\ref#1 #2\par{\setbox0=\hbox{#1}\refindent=\wd0                  %
\plusindent=\refskip                                                 %
\extraindent=\refskip                                                %
\advance\extraindent by 30pt                                         %
\advance\plusindent by -\refindent\tempindent=\parindent %           %
\parindent=0pt\par\hangindent\extraindent %                          %
{#1\hskip\plusindent #2}\parindent=\tempindent}                      %
\refskip=\parindent                                                  %
                                                                     %
%====================================================================%

\def\ol{\overline}

\def\empty{\emptyset}

\def\raj{\restriction}

\def\ol{\overline}

\def\da{\downarrow}

\def\nda{\mathrel{\lower0pt\hbox to 3pt{\kern3pt$\not$\hss}\downarrow}}
\def\nbot{\mathrel{\lower0pt\hbox to 4pt{\kern3pt$\not$\hss}\bot}}
\def\ekom{\mathrel{\lower3pt\hbox to 0pt{\kern3pt$\sim$\hss}\mapsto}}

\def\T{\Theta}
\def\anR{\mathrel{\lower1pt\hbox to 2pt{\kern3pt$R$\hss}\not}}
\def\anoR{\mathrel{\lower1pt\hbox to 2pt{\kern3pt$\overline{R}$\hss}\not}}

\def\anRm{\mathrel{\lower1pt\hbox to 2pt{\kern3pt$R^{-}$\hss}\not}}

\def\ndda{\mathrel{\lower0pt\hbox to 1pt{\kern3pt$\not$\hss}\downdownarrows}}

\null
\vskip 2truecm
\centerline{\bf CONSTRUCTING STRONGLY EQUIVALENT NONISOMORPHIC}
\centerline{\bf MODELS FOR UNSUPERSTABLE THEORIES, PART C}
\vskip 1truecm
\centerline{Tapani Hyttinen and Saharon Shelah$^{*}$}
\vskip 2.5truecm

\chapter{Abstract}

In this paper we prove a strong nonstructure theorem for
$\k(T)$-saturated models of a
stable theory $T$
with dop. This paper continues the work started in [HT].

\vskip 2.5 truecm

\chapter{1. Introduction and basic definitions}

By a strong nonstructure theorem we mean a theorem, which
claims that in a given class of structures, there are
very equivalent nonisomorphic models.
The equivalence is usually measured by the length of Ehrenfeucht-Fraisse
games in which $\exists$ has a winning strategy.
The idea behind this
is, that if models are very equivalent but still nonisomorphic,
they must be very complicated, i.e. there is a lot nonstructure
in the class.

For more background for the theorems of this kind, see [HT].

In this paper we prove the following strong nonstructure theorem
(see Definitions 1.2
and 1.3).

\th 1.1 Theorem. Let $T$ be a stable theory with dop and
$\k =cf(\k )=\l (T)+\k^{<\k (T)}\ge\o_{1}$,
$\l =\l^{<\l}>\k^{+}$ and for all $\xi <\l$,
$\xi^{\k}<\l$. Then there is $F^{a}_{\k}$-saturated model
$M_{0}\models T$ of power $\l$
such that the following is true: for all $\l^{+},\l$-trees $t$
there is a $F_{\k}^{a}$-saturated model
$M_{1}$ of power $\l$ such that
$M_{0}\equiv_{t}^{\l}M_{1}$ and $M_{o}\not\cong M_{1}$.

In [HT] Theorem 1.1 was proved for $F^{a}_{\o}$-saturated
models of a countable superstable theory with dop.
There we used Ehrenfeucht-Mostowski
models to construct
\vskip 1truecm

\noindent
$*$ Research supported by the United States-Israel Binational
Science Foundation. Publ. 602.

\vfill
\eject

\noindent
the required models. To prove that the models
are not isomorphic, it was essential that the sequences in the skeletons
of the models were of finite length. In the case of unsuperstable
theories we cannot quarantee this. Another problem was, of course,
that with Ehrenfeucht-Fraisse models we cannot construct more than
$F^{a}_{\o}$-saturated models.

In this paper we overcome these problems by using $F^{a}_{\k}$-prime models
instead of Ehrenfeucht-Mostowski
models.

\th 1.2 Definition.

(i) Let  $\l$ be a cardinal and $\a$ an ordinal.
Let $t$ be a tree (i.e. for all $x\in t$, the set $\{ y\in t\vert\ y<x\}$
is well-ordered by the ordering of $t$).
If $x,y \in t$ and $\{ z \in t \vert\ z < x \} = \{ z \in t \vert\ z < y \}$,
then we denote $x \sim y$, and the equivalence class
of $x$ for $\sim$ we denote $[x]$.
By a $\l, \a$-tree $t$ we mean a
tree which satisfies:

(a) $\vert [x] \vert < \l$ for every $x \in t$;

(b) there are no branches of length $\ge \a$ in $t$;

(c) $t$ has a unique root;

(d) if $x,y \in t$, $x$ and $y$ have no immediate predecessors
and $x\sim y$, then
$x=y$.

(ii) If $\n$ is a tree and $\a$ is an ordinal then we define the tree
$\a\times\n=(\a\times\n ,<)$
so that $(x,y)<(v,w)$ iff $y<w$ or $y=w$ and $x<v$.

\th 1.3 Definition. Let $t$ be a tree and $\k$ a cardinal.
The  Ehrenfeucht-Fraisse game
of length $t$ between models $\A$ and $\B$,
$G^{\k}_{t}(\A, \B)$, is the following.
At each move $\a$:

(i) player $\forall$ chooses $x_\a \in t$, $\k_{\a}<\k$ and
either $a_\a^\b \in \A$, $\b <\k_{\a}$ or $b_\a^\b \in \B$,
$\b <\k_\a$, we will denote this sequence
by $X_{\a}$;

(ii) if $\forall$ chose from $\A$ then
$\exists$ chooses $b_\a^\b \in \B$, $\b <\k_\a$, else
$\exists$ chooses
$a_\a^\b \in \A$, $\b <\k_\a$, we will denote this sequence by $Y_{\a}$.

\noindent
$\forall$ must move so that $(x_\b)_{\b \le \a}$
form a strictly increasing sequence in $t$.
$\exists$ must move so that
$\{ (a_\g^\b, b_\g^\b) \vert \g \le \a , \b <\k_\g \}$
is a partial isomorphism from $\A$ to $\B$.
The player who first has to break the rules loses.

We write $\A\equiv^{\k}_{t}\B$ if $\exists$ has a winning strategy
for $G^{\k}_{t}(\A ,\B )$.

The following theorem is frequently used in this paper.

\th 1.4 Theorem. ([Sh]) Let $T$ be a stable theory.
Assume $I$ is an infinite indiscernible
sequence over $A$, $I\subseteq B$ and $J\subseteq I$
is countable.

(i) $Av(I,B)$ does not fork over $J$ and $Av(I,J)$ is stationary.

(ii) $I\cup\{ a\}$ is indiscernible over $A$ iff
$t(a,A\cup I)=Av(I,A\cup I)$.

\proof See [Sh] Lemma III 4.17. $\eop$

\th 1.5 Corolary. Let $T$ be a stable theory. Assume $I$ is an infinite
indiscernible sequence over $A$ and $J\subseteq I$ is infinite.
Then $I-J$ is independent over $A\cup J$.

\proof Follows immediately from Theorem 1.4. $\eop$

\chapter{2. Construction}

Through out this paper we assume that
$T$ is a stable theory with dop,
$\k =cf(\k )=\l (T)+\k^{<\k (T)}\ge\o_{1}$,
$\l =\l^{<\l}>\k^{+}$ and for all $\xi <\l$,
$\xi^{\k}<\l$.

\th 2.1 Theorem. ([Sh]) There are models $\A_{i}$, $i<3$,
of cardinality $<\k$ and infinite indiscernible
sequence $I$ over $\A_{1}\cup\A_{2}$ such that

(i) $\A_{0}\subseteq\A_{1}\cap\A_{2}$, $\A_{1}\da_{\A_{0}}\A_{2}$,

(ii) $Av(I,I\cup\A_{1}\cup\A_{2})\perp\A_{1}$,
$Av(I,I\cup\A_{1}\cup\A_{2})\perp\A_{2}$,

(iii) $t(I,\A_{1}\cup\A_{2})$ is almost orthogonal to $\A_{1}$
and to $\A_{2}$,

(iv) if $B_{i}$, $i<3$ are such that $B_{0}\da_{\A_{0}}\A_{1}\cup\A_{2}$,
$B_{1}\da_{\A_{1}\cup B_{0}}\A_{2}\cup B_{2}$
and $B_{2}\da_{\A_{3}\cup B_{0}}\A_{1}\cup B_{1}$
then
$$t(I,\A_{1}\cup\A_{2})\vdash t(I,\A_{1}\cup\A_{2}\cup\bigcup_{i<3}B_{3}).$$

\proof This is [Sh] X Lemma 2.4, except that in (iv), only
$$(*)\ \ \ \ stp(I,\A_{1}\cup\A_{2})
\vdash t(I,\A_{1}\cup\A_{2}\cup\bigcup_{i<3}B_{3})$$
is proved. But since $\k\ge\k_{r}(T)$, by [Sh] XI Lemma 3.1
$\A_{1}\cup\A_{2}$ is a good set. It is easy to see that
this together with (*) implies
$$t(I,\A_{1}\cup\A_{2})\vdash t(I,\A_{1}\cup\A_{2}\cup\bigcup_{i<3}B_{3}).$$
$\eop$

In [HT] the following theorem is proved.

\th 2.2 Theorem. ([HT] Theorem 3.4)
There is a $\l^{+},\l +1$-tree $\n$ such that
it has a branch of length $\l$ and
for every $\l^{+},\l$-tree $t$ there is a $\l^{+},\l$-tree $\xi$
such that $\n\equiv_{t}^{\l}\xi$.

Let $\n$ be a tree. We define a model $M(\n)$. Let $\A ,\B ,\C$
and $I$ be as $\A_{0}, \A_{1}, \A_{2}$ and $I$
in Theorem 2.1. We may assume that $\vert I\vert =\l$.

For all $t\in\n$ we choose $\A_{t},\B_{t}$ and $\C_{t}$ so that

(i) there is an automorphism $f_{t}$ (of the monster model) such
that $f_{t}(\B_{t})=\B$, $f_{t}(\C_{t})=\C$ and $f^{-1}_{t}\raj\A =id_{\A}$,

(ii) $\B_{t}\cup\C_{t}\da_{\A}
\bigcup\{\B_{s}\cup\C_{s}\vert\ s\in\n ,\ s\ne t\}$.

\noindent
For all $s,t\in\n$, $s<t$, we choose $I_{st}$ so that

(i) there is an automorphism $g_{st}$ such that
$g_{st}\raj\B_{s}=f_{s}\raj\B_{s}$, $g_{st}\raj\C_{t}=f_{t}\raj\C_{t}$
and $g_{st}(I_{st})=I$,

(ii) $I_{st}\da_{\B_{s}\cup\C_{t}}
\bigcup\{\B_{p}\cup\C_{p}\vert
\ p\in\n\}\cup\bigcup\{ I_{pr}\vert
\ p,r\in\n ,\ p<r,\ p\ne s\ \hbox{\rm or}\ r\ne t\}$.

We define $M(\n )$ to be the $F^{a}_{\k}$-primary model
over $S(\n )=\bigcup\{\B_{t}\cup\C_{t}\vert\ t\in\n\}\cup\bigcup
\{ I_{st}\vert\ s,t\in\n ,\ s<t\}$.

By Theorem 2.2,
Theorem 1.1 follows immediately from the theorem below.

\th 2.3 Theorem. Let $\n$ be as in Theorem 2.2 and
$M_{0}=M(\n )$. Assume $t$ is a $\l^{+},\l$-tree.
Let $\xi$ be a $\l^{+},\l$-tree such that
$\n\equiv_{\k\times t}^{\l}\xi$.
If $M_{1}=M(\xi )$,
then $M_{0}\equiv_{t}^{\l}M_{1}$, $M_{o}\not\cong M_{1}$
and the cardinality of the models is $\l$.

The claim on the cardinality of the models follows immediately
from the assumptions on $\l$. The other two claims are
proved in the next two chapters.

Notice that in $\xi$ there are no brances of length $\l$. Since in
$\n$ there is such a branch, this enables us to prove the nonisomorphism
of the models.

\chapter{3. Equivalence}

In this chapter we prove the first part of Theorem 2.3. We want to remind
the reader of the assumptions made in the beginning of Chapter 2.

Let $(S(\n ),\{ d_{i}\vert\ i<\a\},(D_{i}\vert\ i<\a ))$
and
$(S(\xi ),\{ e_{i}\vert\ i<\a\},(E_{i}\vert\ i<\b ))$
be $F^{a}_{\k}$-constructions of $M(\n )$ and $M(\xi )$,
respectively, see [Sh] IV Definition 1.2.
If we choose the constructions
carefully we can assume $\a =\b =\l$.

We enumerate $\n$ and $\xi$: $\n =\{ t^{\n}_{i}\vert\ i<\l\}$
and $\xi =\{ t^{\xi}_{i}\vert\ i<\l\}$. Furthermore we do this
so that if $t^{*}_{i}<t^{*}_{j}$ then $i<j$, $*\in\{\n ,\xi\}$.
If $\g\le\l$, we write $\n (\g )=\{ t^{\n}_{i}\vert\ i<\g\}$
and similarly for $\xi (\g )$.

We also enumerate all $I_{st}$: $I_{st}=\{ a_{st}^{i}\vert\ i<\l\}$.

We write $S(\n ,\g )$ for
$$\bigcup\{\B_{t}\vert\ t\in \n (\g )\}
\cup\bigcup\{\C_{t}\vert\ t\in \n (\g )\}\cup$$
$$\bigcup\{ a_{st}^{i}\vert\ s<t,\ s,t\in\n (\g ),\ i<\g\}$$
and similarly for $S(\xi ,\g )$.

If $\g <\l$ and $g:\n (\g )\rightarrow\xi (\g )$ is
a partial isomorphism then by $g^{*}$ we mean the function
from $S(\n ,\g )$ onto $S(\xi ,\g )$ which satisfies:

(i) if $g(t)=t'$ then for all $a\in B_{t}$ and $b\in C_{t}$,
$g^{*}(a)=f^{-1}_{t'}(f_{t}(a))$ and $g^{*}(b)=f^{-1}_{t'}(f_{t}(b))$,

(ii) if $g(t)=t'$, $g(s)=s'$, $t<s$ and $a\in I_{ts}$ then
$g^{*}(a)=g^{-1}_{t's'}(g_{ts}(a))$.

\th 3.1 Lemma. If $\g <\l$ and $g:\n (\g )\rightarrow\xi (\g )$ is
a partial isomorphism then $g^{*}$ is a partial isomorphism.

\proof Immediate by the definitions. $\eop$

We write
$$M(\n ,\g )=S(\n ,\g )\cup\{ d_{i}\vert\ i<\g\}$$
and
similarly for $M(\xi ,\g )$.
We say that $\g <\l$ is good if for all $i<\g$, $D_{i}\subseteq M(\n ,\g )$
and $E_{i}\subseteq M(\xi ,\g )$. Notice that the set of
all good ordinals is cub in $\l$. Notice also that the set of
those ordinals $\g <\l$ for which $M(\n ,\g )$ is
$F^{a}_{\k}$-saturated, is $\ge\k$-cub, i.e. it is unbounded
in $\l$ and closed under increasing sequences of cofinality $\ge\k$.

\th 3.2 Lemma. Assume $A\subseteq B$, $a_{i}$ and $C_{i}$, $i<\a$,
are such that

(i) $C_{i}\subseteq A\cup\{ a_{j}\vert\ j<i\}$ is of power $<\k$,

(ii) $t(a_{i},C_{i})\vdash t(a_{i},B\cup \{ a_{j}\vert\ j<i\})$.

\noindent
Then for all sequences $\ol d\in\{ a_{i}\vert\ i<\a\}$, there is
$D\subseteq A$ of power $<\k$ such that
$t(\ol d ,D)\vdash t(\ol d ,B)$.
Especially, $\ol d\da_{A}B$.

\proof See the proof of [Sh] Theorem IV 3.2. $\eop$

\th 3.3 Lemma. Let $\g <\l$ be good, $\g <\d <\l$,
$g:\n (\d )\rightarrow\xi (\d )$
is a partial isomorphism, $f:M(\n ,\g )\rightarrow M(\xi ,\g )$
is a partial isomorphism and $g^{*}\raj S(\n ,\g )\subseteq f$. Then
$f\cup g^{*}$ is a partial isomorphism from $M(\n ,\g )\cup S(\n ,\d )$
onto $M(\xi ,\g )\cup S(\xi ,\d )$.

\proof Follows immediately from Lemmas 3.1, 3.2 and the definition of a good
ordinal. $\eop$

\th 3.4 Lemma. Assume $\g <\l$ is good, $g:\n (\g )\rightarrow\xi (\g )$
and $f:M(\n ,\g )\rightarrow M(\xi ,\g )$ are partial isomorphism,
$g^{*}\subseteq f$ and
$$(\n ,a)_{a\in\n (\g )}\equiv^{\l}_{\k}(\xi ,f(a))_{a\in \n (\g )}.$$
If $A\subseteq M_{0}$ is of power $<\l$ then there are
good $\g' <\l$, partial isomorphisms $g':\n (\g' )\rightarrow\xi (\g' )$
and $f':M(\n ,\g' )\rightarrow M(\xi ,\g' )$ such that
$(g')^{*}\subseteq f'$, $f\subseteq f'$, $g\subseteq g'$ and
$A\subseteq M(\n ,\g' )$.

\proof By playing the Ehrenfeucht-Fraisse game we can find
a good $\g'<\l$ such that

(i) there is a partial isomorphism $g':\n (\g')\rightarrow\xi (\g')$
such that $g\subseteq g'$,

(ii) $M(\n ,\g')$ is $F^{a}_{\k}$-primary over $S(\n ,\g')$ and
$M(\xi ,\g')$ is $F^{a}_{\k}$-primary over $S(\xi ,\g')$,

(iii) $A\subseteq M(\n ,\g')$.

By (i) above and Lemma 3.3,
$f\cup (g')^{*}$ is a partial isomorphism from $M(\n ,\g )\cup S(\n ,\g' )$
onto $M(\xi ,\g )\cup S(\xi ,\g')$.
\relax From (ii) it follows that
$M(\n ,\g')$ is $F^{a}_{\k}$-primary over $M(\n ,\g )\cup S(\n ,\g')$ and
$M(\xi ,\g')$ is $F^{a}_{\k}$-primary over $M(\xi ,\g )\cup S(\xi ,\g')$.
So the existence of the required $f'$ follows from the uniqueness of
the $F^{a}_{\k}$-primary models ([Sh] Conclusion IV 3.9). $\eop$

\th 3.5 Theorem. $M_{0}\equiv^{\l}_{t}M_{1}$.

\proof By Lemma 3.4, it is easy to translate the winning strategy
of $\exists$ in $G^{\l}_{\k\times t}(\n ,\xi )$ to her winning strategy
in $G^{\l}_{t}(M_{0},M_{1})$.
$\eop$

\chapter{4. Nonisomorphism}

In this chapter we prove the second part of Theorem 2.3, i.e.
$M_{0}\not\cong M_{1}$. Again we want to remind the reader
of the assumptions made in the beginning of Chapter 2.

For a contradiction we assume that $f:M_{0}\rightarrow M_{1}$
is an isomorphism.

If $a\in M_{0}$ then we write $\a_{a}$ for
the least $\a$ such that $a\in M(\n ,\a )$ and similarly for $a\in M_{1}$.
By $\a_{A}$ we mean $\bigcup\{ \a_{a}\vert\ a\in A\}$.

Let $X\subseteq\n$ be such that $\vert X\vert =\l$ and
for all $x,y\in X$ if $x\ne y$ then either $x<y$ or
$y<x$.
For every $x\in X$
we choose $u_{x}^{i}$, $S_{x}^{i}$ and $N_{x}^{i}$, $i\in\{ 0,1\}$, so that

(i) $x\in u_{x}^{0}\subseteq\n$ and $u_{x}^{1}\subseteq\xi$,

(ii) $S_{x}^{i}=\bigcup\{\B_{t}\vert\ t\in u_{x}^{i}\}\cup
\bigcup\{\C_{t}\vert\ t\in u_{x}^{i}\}\cup
\bigcup\{ I_{st}^{x}\vert\ s,t\in u_{x}^{i},\ s<t\}$,
where $I_{st}^{x}\subseteq I_{st}$ is of infinite power at most
$\k$,

(iii) $N_{x}^{i}\subseteq M_{i}$ is $F^{a}_{\k}$-primary over $S_{x}^{i}$
and furthermore if $a\in N_{x}^{0}-S(\n )$ and $a=d_{i}$ in the
construction of $M_{0}$ then $D_{i}\subseteq N_{x}^{0}$ and
similarly for $N_{x}^{1}$,

(iv) $f\raj N_{x}^{0}$ is onto $N_{x}^{1}$,

(v) $\vert N_{x}^{i}\vert\le\k$,

(vi) if $M(\n ,\a )$ is $F^{a}_{\k}$-saturated, then so is
$M(\n ,\a )\cap N_{x}^{0}$.

\noindent
It is easy to see that these sets exist.

\th 4.1 Lemma. Assume $A_{i}$, $i<\l$, are sets of power $\le\k$.
Then there are $X\subseteq\l$ and $B$ such that $\vert X\vert =\l$
and for all $i,j\in X$, $A_{i}\cap A_{j}=B$.

\proof Without loss of generality we may assume that for all
$i<\l$, $A_{i}\subseteq\l$. We define
$f(\a )=sup(A_{i}\cap(\cup_{j<i}A_{j}))$.
Since $\l >\k^{+}$ is regular,
this function is regressive on a stationary set.
So by Fodor's lemma,
it is constant on some set $X'$ of power $\l$. Since for all
$\theta <\l$, $\theta^{\k}<\l$, the claim follows by the pigeon hole
principle. $\eop$

By Lemma 4.1 and the pigeon hole principle
we may assume that $X$ is chosen so that it satisfies
the following:

(i) There are $u^{i}$, $S^{i}$ and $N^{i}$, $i\in\{ 1,2\}$,
such that for all $x,y\in X$,
if $x\ne y$ then $u^{i}_{x}\cap u^{i}_{y}=u^{i}$,
$S^{i}_{x}\cap S^{i}_{y}=S^{i}$ and
$N^{i}_{x}\cap N^{i}_{y}=N^{i}$.

(ii) For all $x\in X$,
$M(\n ,\a_{N^{0}})\cap N^{o}_{x}=N^{0}$
and if $x<y$ then
$M(\n ,\a_{N^{0}_{x}})\cap N^{o}_{y}=N^{0}$ and
similarly for $1$ instead of $0$.

(iii) For all $x,y\in X$, there are elementary maps
$f^{i}_{xy}:N^{i}_{x}\rightarrow N^{i}_{y}$ and
an order isomorphisms
$g^{i}_{xy}:u^{i}_{x}\rightarrow u^{i}_{y}$ such that

(a) $f^{i}_{xy}\raj N^{i}=id_{N^{i}}$, $g^{i}_{xy}\raj u^{i}=id_{u^{i}}$
and $g^{0}_{xy}(x)=y$,

(b) for all $t\in u^{i}_{x}$ and $a\in \B_{t}\cup\C_{t}$,
$f^{i}_{xy}(a)=f^{-1}_{g^{i}_{xy}(t)}(f_{t}(a))$,

(c) for all $s,t\in u^{i}_{x}$, $s<t$,
$f^{i}_{xy}\raj I^{x}_{st}$ is onto $I^{y}_{g^{i}_{xy}(s)g^{i}_{xy}(t)}$

(d) for all $a\in N^{0}_{x}$, $f(f^{0}_{xy}(a))=f^{1}_{xy}(f(a))$.

\th 4.2 Lemma. Let $x,y\in X$, $x<y$.

(i) $N^{i}$ is $F^{a}_{\k}$-primary over $S^{i}$.

(ii) $N^{i}_{x}$ is $F^{a}_{\k}$-primary over $N^{i}\cup S^{i}_{x}$.

(iii) $N^{i}\da_{S^{i}}S^{i}_{x}\cup S^{i}_{y}$.

(iv) $N^{i}_{x}\da_{N^{i}}N^{i}_{y}$.

(v) $I_{xy}\da_{\B_{x}\cup\C_{y}}N^{o}_{x}\cup N^{o}_{y}$.

\proof Immediate by (ii) in the choice of $X$ and Lemma 3.2. $\eop$

\th 4.3 Corollary. Let $x,y\in X$, $x<y$.

(i) If $A,B$ and $C$ are such that $A\da_{N^{0}}N^{0}_{x}\cup N^{0}_{y}$,
and $B\cup N^{0}_{x}\da_{N^{0}\cup A}N^{0}_{y}\cup C$ then
$$t(I_{xy},\B_{x}\cup\C_{y})\vdash
t(I_{xy},I_{xy}\cup N^{o}_{x}\cup N^{0}_{y}\cup A\cup B\cup C)$$

(ii) $t(I_{xy}\cup N^{0}_{x}\cup N^{0}_{y},\empty )$ does not depend on
$x$ and $y$.

\proof (i) By the first assumption on $A$ and Lemma 4.2 (iii)
$$A\cup N^{0}\da_{S^{0}}\B_{x}\cup\C_{y}.$$
By the construction of $\M_{0}$, this implies
$$(A)\ \ \ \ \ A\cup N^{0}\da_{\A}\B_{x}\cup\C_{y}.$$
\relax From the second assumption it follows easily that
$$(B)\ \ \ \ \ B\cup N^{0}_{x}\da_{\B_{x}\cup N^{0}\cup A}N^{0}_{y}\cup C$$
and
$$(C)\ \ \ \ \ C\cup N^{0}_{y}\da_{\C_{y}\cup N^{0}\cup A}N^{0}_{x}\cup B.$$
By Theorem 2.1 (iv), (A),(B) and (C) imply the claim.

(ii) By (iii) in the choice of $X$ and Lemma 4.2 (iv),
for all $x'<y'$, $f^{0}_{x x'}\cup f^{0}_{y y'}$ is an elementary
map. So the claim follows from (A), (B) and (C) above and
Theorem 2.1 (iv).
$\eop$

For $x,y\in X$, $x<y$,
let $I^{c}_{xy}$ be some countable subset of $I_{xy}$.

\th 4.4 Lemma. Assume $x,y\in X$, $x<y$. Then there are
$s\in u^{1}_{x}-u^{1}$
and $t\in u^{1}_{y}-u^{1}$ such that either

(i) $s<t$ and
$Av(f(I^{c}_{xy}),f(I^{c}_{xy}\cup \B_{x}\cup \C_{y}))$
is not orthogonal to
$Av(I^{c}_{st},I^{c}_{st}\cup\B_{s}\cup\C_{t})$,

\noindent
or

(ii) $t<s$ and
$Av(f(I^{c}_{xy}),f(I^{c}_{xy}\cup \B_{x}\cup \C_{y}))$
is not orthogonal to
$Av(I^{c}_{ts},I^{c}_{ts}\cup\B_{t}\cup\C_{s})$.

\proof For a contradiction, we assume that such $s$ and $t$ do not
exist.

Let

$\xi^{0} (x,y)=\{ (s,t)\vert\ s<t\ \hbox{and}
\ s\in u^{1}_{x}-u^{1},\ t\not\in u^{1}_{y}-u^{1}
\ \hbox{or}\ t\in u^{1}_{x}-u^{1},\ s\not\in u^{1}_{y}-u^{1}\}$

$\xi^{1} (x,y)=\{ (s,t)\vert\ s<t\ \hbox{and}
\ s\not\in u^{1}_{x}-u^{1},\ t\in u^{1}_{y}-u^{1}
\ \hbox{or}\ t\not\in u^{1}_{x}-u^{1},\ s\in u^{1}_{y}-u^{1}\}$ and

$\xi^{2} (x,y)=\{ (s,t)\vert\ s<t\ \hbox{and}
\ s\in u^{1}_{x}-u^{1},\ t\in u^{1}_{y}-u^{1}
\ \hbox{or}\ t\in u^{1}_{x}-u^{1},\ s\in u^{1}_{y}-u^{1}\}$.
For $i\in\{ 0,1,2\}$, let

$S^{i}(x,y)=S(\xi )-
(S^{1}_{x}\cup S^{1}_{y}\cup
\bigcup_{j\ge i}\{ I_{st}\vert\ (s,t)\in\xi^{j}(x,y)\} )$

\noindent
and

$R^{i}(x,y)=
\{ I_{st}\vert\ (s,t)\in\xi^{i}(x,y)\}$.

Now it is easy to see that
$S^{0}(x,y)\da_{S^{1}}S^{1}_{x}\cup S^{1}_{y}$.
By Lemma 3.2 $N^{1}\da_{S^{1}}
S^{0}(x,y)\cup S^{1}_{x}\cup S^{1}_{y}$. So
$$S^{0}(x,y)\da_{N^{1}}S^{1}_{x}\cup S^{1}_{y}.$$
By Lemma 4.2 this implies
$$(A)\ \ \ S^{0}(x,y)\da_{N^{1}}N^{1}_{x}\cup N^{1}_{y}.$$

By the construction
$$(B)\ \ \ R^{0}(x,y)\cup S^{1}_{x}\da_{S^{1}\cup S^{0}(x,y)}
R^{1}(x,y)\cup S^{1}_{y}.$$
By Lemma 3.2
$$N^{1}_{x}\da_{S^{1}_{x}}S^{0}(x,y)\cup R^{0}(x,y)
\cup R^{1}(x,y)\cup S^{1}_{y}$$
and so
$$R^{0}(x,y)\cup N^{1}_{x}\da_{S^{1}_{x}\cup S^{0}(x,y)\cup R^{0}(x,y)}
R^{1}(x,y)\cup S^{1}_{y}.$$
By (B) this implies
$$(C)\ \ \ R^{0}(x,y)\cup N^{1}_{x}\da_{N^{1}\cup S^{0}(x,y)}
R^{1}(x,y)\cup S^{1}_{y}.$$
By Lemma 3.2 and (ii) in the choice of $X$,
$$N^{1}_{y}\da_{S^{1}_{y}}S^{0}(x,y)\cup R^{0}(x,y)
\cup R^{1}(x,y)\cup N^{1}_{x}$$
and so
$$R^{1}(x,y)\cup N^{1}_{y}\da_{S^{1}_{y}\cup S^{0}(x,y)\cup R^{1}(x,y)}
R^{0}(x,y)\cup N^{1}_{x}.$$
By (C) this implies
$$(D)\ \ \ R^{1}(x,y)\cup N^{1}_{y}\da_{N^{1}\cup S^{0}(x,y)}
R^{0}(x,y)\cup N^{1}_{x}.$$

Then by (A), (D) and Corollary 4.3 (i), $f(I_{xy})$ is indiscernible over
$N^{1}_{x}\cup N^{1}_{y}\cup S^{2}(x,y)$.

By Lemma 3.2 and (ii) in the choice of $X$, we see that
for all $(s,t)\in\xi^{2} (x,y)$,
$I_{st}$ is indiscernible over $N^{1}_{x}\cup N^{1}_{y}\cup S^{2}(x,y)$
and $(I_{st})_{(s,t)\in\xi^{2}(x,y)}$ is independent over
$N^{1}_{x}\cup N^{1}_{y}\cup S^{2}(x,y)$.

For all $(u,v)\in\xi^{2}(x,y)\cup\{ (x,y)\}$
we choose infinite $I^{*}_{uv}\subseteq I_{uv}$
of power $<\l$ such that

(i) for all $(u,v)\in\xi^{2}(x,y)$, if we write
$B(u,v)=N^{1}_{x}\cup N^{1}_{y}\cup S^{2}(x,y)\cup I^{*}_{uv}$,
then
$$I_{uv}-I_{uv}^{*}\da_{B(u,v)}f(I^{*}_{xy})
\cup\bigcup\{ I^{*}_{st}\vert\ (s,t)\in\xi^{2}(x,y),\ (s,t)\ne (u,v)\} .$$

(ii) $I^{c}_{xy}\subseteq I^{*}_{xy}$ and if we write
$B(x,y)=N^{1}_{x}\cup N^{1}_{y}\cup S^{2}(x,y)\cup f(I^{*}_{xy})$,
then
$$f(I_{xy}-I_{xy}^{*})\da_{B(x,y)}
\bigcup\{ I^{*}_{st}\vert\ (s,t)\in\xi^{2}(x,y)\} .$$

\noindent
Because $\vert \xi^{2}(x,y)\vert <\l$, it is easy to see that such
$I^{*}_{uv}$ exist.

Since $Av(f(I^{c}_{xy}),f(I^{c}_{xy}\cup \B_{x}\cup \C_{y}))$
is orthogonal to
$Av(I^{c}_{st},I^{c}_{st}\cup\B_{s}\cup\C_{s})$ for all $(s,t)\in\xi^{2} (x,y)$
we see that $I_{xy}-I^{*}_{xy}$ is indiscernible over $S(\xi )$.
Because $\vert I_{xy}-I^{*}_{xy}\vert =\l$,
this contradicts [Sh] Theorem IV 4.9 (2). $\eop$

If $s,t\in\xi$, then we write $\T_{st}$ for the set of all
infinite $J$ such that for some $J'$,
$J\subseteq J'$ and
there is an automorphism
$g$ for which $g\raj\B_{s}=f_{s}\raj\B_{s}$, $g\raj\C_{t}=f_{t}\raj\C_{t}$
and $g(J')=I$.

\th 4.5 Lemma. Assume $x,y\in X$, $x<y$, $s\in u^{1}_{x}-u^{1}$,
$t\in u^{1}_{y}-u^{1}$ and $s$ and $t$ are incomparable in $\xi$.
If $J\in\T_{st}$,
then $Av(f(I^{c}_{xy}),f(I^{c}_{xy}\cup \B_{x}\cup \C_{y}))$ is orthogonal to
$Av(J,J\cup\B_{s}\cup\C_{t})$. Also if $J\in\T_{ts}$,
then $Av(f(I^{c}_{xy}),f(I^{c}_{xy}\cup \B_{x}\cup \C_{y}))$ is orthogonal to
$Av(J,J\cup\B_{t}\cup\C_{s})$.

\proof For a contradiction assume that
$Av(f(I^{c}_{xy}),f(I^{c}_{xy}\cup \B_{x}\cup \C_{y}))$ is not orthogonal to
$Av(J,J\cup\B_{s}\cup\C_{t})$, the other case is similar.
Then we can
choose $J$ so that in addition, $\vert J\vert =\o$ and
$J\subseteq M_{1}$.

By Theorem 2.1 (iv), $J$ is indiscernible over
$S(\xi )$. By [Sh] Theorem IV 4.14,
$Av(J,M_{1})$ is $F^{a}_{\k^{+}}$-isolated.
Then we can find
a model $D\subseteq M_{1}$ of power $\le\k$ such that

(a) $f(I_{xy}^{c}\cup\B_{x}\cup\C_{y})\cup J\cup\B_{s}\cup\C_{t}\subseteq D$,

(b) $Av(f(I^{c}_{xy}),D)$ is not almost orthogonal to $Av(J,D)$,

(c) $Av(J,D)\vdash Av(J,M_{1})$.

\noindent
(For (c), notice that because $D$ is a model,
$t(a,D)\vdash stp(a,D)$.)
But since $\vert D\vert <\l$ and $\vert f(I_{xy})\vert =\l$,
it is easy to see that $Av(f(I^{c}_{xy}),D)$
is satisfied in $M_{1}$, a contradiction. $\eop$

Let $x,y\in X$ be such that $x<y$. By Lemma 4.4
we can find $s_{xy}$ and $t_{xy}$
such that there is $J\in\T_{s_{xy}t_{xy}}\cup\T_{t_{xy}s_{xy}}$
for which $Av(f(I^{c}_{xy}),f(I^{c}_{xy}\cup \B_{x}\cup \C_{y}))$
is not orthogonal to
$Av(J,J\cup\B_{s_{xy}}\cup\C_{t_{xy}})$ or to
$Av(J,J\cup\B_{t_{xy}}\cup\C_{s_{xy}})$.
By Lemma 4.3 (ii) we can choose these
so that for all $y$ and $y'$ from $X$,
if $x<y$ and $x<y'$ then $s_{xy}=s_{xy'}$.
We call this element just $s_{x}$.
Similarly we can choose $t_{xy}$ so that it does not depend on $x$ ($x<y$).
We call this element $t_{y}$.

\th 4.6 Lemma. For all $x$ and $x'$ from $X$,
$s_{x}$ and $s_{x'}$ are comparable in $\xi$.

\proof By Lemma 4.5, for all $y\in X$, if $y>x$ and $y>x'$ then
$t_{y}$ is comparable to $s_{x}$ and to $s_{x'}$.
Since $\vert\{ z\in \xi\vert\ z\le s_{x}\ \vee\ z\le s_{x'}\}\vert <\l$
and if $y\ne y'$ then $t_{y}\ne t_{y'}$,
we can find $y\in X$ such that
$s_{x}<t_{y}$ and $s_{x'}<t_{y}$, which implies the claim. $\eop$

\th 4.7 Theorem. $M_{0}\not\cong M_{1}$.

\proof If $M_{0}\cong M_{1}$ then by Lemma 4.6 we can find
$Y\subseteq\xi$ of power $\l$ such that for all $s,t\in Y$
if $s\ne t$ then
either $s<t$ or $t<s$. Clearly this contradicts the fact that $\xi$
is a $\l^{+},\l$-tree.
$\eop$

Together with Theorem 3.5, Theorem 4.7 implies Theorem 2.3,
and so Theorem 1.1 is proved.

\th 4.8 Remark. As in [HT], we can see that Theorem 1.1 implies
the following: Under the assumptions of Theorem 1.1, for every
$\l^{+},\l$-tree $t$ there are models $M_{i}\models T$, $i<\l^{+}$,
such that for all $i<j<\l^{+}$, $M_{i}\equiv^{\l}_{t}M_{j}$
and $M_{i}\not\cong M_{j}$.

\chapter{References.}

\item{[HT]} T. Hyttinen and H. Tuuri, Constructing strongly equivalent
nonisomorphic models for unstable theories, APAL 52, 1991, 203-248.

\item{[Sh]} S. Shelah, Classification Theory, Stud. Logic Found. Math. 92
(North-Holland, Amsterdam, 2nd rev. ed., 1990).

\medskip

Tapani Hyttinen

Department of Mathematics

P.O. Box 4

00014 University of Helsinki

Finland

\medskip

Saharon Shelah

Institute of Mathematics

The Hebrew University

Jerusalem

Israel

\medskip

Rutgers University

Hill Ctr-Bush

New Brunswick

New Jersey 08903

U.S.A.

\end